\documentclass[a4paper, 10pt]{article}
\usepackage{amsmath,amssymb,amsthm}
\usepackage{enumerate}
\usepackage{graphicx}

\newcommand{\R}{\mathbb{R}}
\newcommand{\N}{\mathbb{N}}
\newcommand{\C}{\mathbb{C}}
\newcommand{\Z}{\mathbb{Z}}

\newcommand{\G}{\mathbb{G}}

\def\ben{\begin{enumerate}}
\def\een{\end{enumerate}}

\def\ga{\,^\gamma\hskip-1pt} 
\def\sa{\,^s\hskip-1pt}

\def\demo{{\bf Proof}: }

\newtheorem{theorem}{Theorem}[section]
\newtheorem{proposition}[theorem]{Proposition}
\newtheorem{lemma}[theorem]{Lemma}

\newtheorem{definition}[theorem]{Definition}

\begin{document}

\title{Infinite families of inequivalent real circle actions on affine four-space}
\author{L. Moser-Jauslin
}
\date{}

\maketitle

\abstract{The main result of this article is to construct infinite families of non-equivalent equivariant real forms of linear $\C^*$-actions on affine four-space. We consider the real form of $\C^*$ whose fixed point is a circle. In \cite{F-MJ} one example of a non-linearizable circle action was constructed. Here,  this result is generalized by developing a new approach which allows us to compare different real forms. The constructions of these forms are based on  the structure of equivariant $O_2(\C)$-vector bundles.}

\section{Introduction}

In this article, we construct  infinite families of inequivalent real forms of linear $\C^*$-actions on complex affine four-space. Consider the real form $\sigma$ of $\C^*$ whose fixed point set is a circle.  That is, $\sigma(t)=\overline{t}^{-1}$ for all $t\in\C^*$. Now consider the linear action of $\C^*$  on the four dimensional affine space $Y$ with weights $2,-2,n,-n$ for $n$ an odd integer $\ge 3$. We study equivariant real forms of this $\C^*$-variety, that is, real forms  $\mu$ of $Y$ which are compatible with the real form $\sigma$ of $\C^*$.  The construction of these real forms uses a classification of equivariant $O_2(\C)$-vector bundles where the base (resp. zero fiber) is the irreducible $O_2(\C)$-module for which $\C^*$ acts with weights $2$ and $-2$ (resp. $n$ and $-n$). For $n\ge 5$, it is known that one can find a family for which the actions of $O_2(\C)$ on the total spaces are all inequivalent. By adapting these ideas to real forms, we will show in particular that, for $n\ge 5$, there are infinitely many distinct real forms of $Y$.

It was shown in \cite{K-R} that, over any field of characteristic zero  all actions of any form of a $\G_m$-action  on affine three-space is linearizable.  The result from \cite{F-MJ} shows that this result does not hold in dimension 4 for the field $\R$. Here,  we describe a systematic approach to find many non-linearizable actions which are all pairwise inequivalent.

In   section  2, we recall basic definitions of real forms of equivariant varieties. Then in section 3, we state and give the proof of the main results. In section 4 and 5, we describe the method of constructing the examples, using equivariant vector bundles for the group $O_2(\C)$. We then show in section 6  how the analogous result for another case fails to hold. In the final section, we state two open questions related to these results.

\section{Description of real forms}

If $G$ is an algebraic linear complex group, a real form of $G$ is given by an antiholomorphic involution $\sigma$ on $G$ where $\sigma$ is a group automorphism.

Given a $G$-variety $Y$, a $(G,\sigma)$-real form of $Y$, or an \emph{equivariant real form of $Y$ compatible with $\sigma$},   is given by an antiholomorphic involution $\mu$ of $Y$ such that $\mu(gy)=\sigma(g)\mu(y)$ for all $g\in G$ and $y\in Y$. Two real forms $\mu_1$ and $\mu_2$ are equivalent if there exists a regular $G$-automorphism $\psi$ of $Y$ such that $\mu_2=\psi\circ\mu_1\circ\psi^{-1}$. 

Note that, if one such real form exists, then the set of all equivalence classes of real forms is determined by a cohomology set. More precisely, let $\Gamma=\{1,\gamma\}$ be the Galois group of $\C/\R$. Let $\mu_0$ be one real $(G,\sigma)$-real form of the $G$-variety $Y$. Consider the action of $\Gamma$ on $A=Aut_G(Y)$ given by : $\ga\psi=\mu_0\circ\psi\circ\mu_0$. The set of real $(G,\sigma)$-forms is determined by $Z^1(\Gamma,A)=\{\varphi\in A: \ga\varphi=\varphi^{-1}\}$. More precisely, the set of real $(G,\sigma)$-forms are given by $\varphi\circ\mu_0$, where $\varphi\in Z^1(\Gamma,A)$. 

 Consider the equivalence relation defined by : $\varphi_1\sim\varphi_2$ if  there exists $\psi\in A$ such that $\varphi_2=\psi\circ\varphi_1\circ\ga (\psi)^{-1}$. 
 Thus, $H^1(\Gamma,A)=Z^1(\Gamma,A)/\sim$ is in bijection with the equivalence classes of $(G,\sigma)$-real forms of $Y$.

\section{Main result : the $(2,2m+1)$ case}

\subsection{Construction of circle forms}

The case considered here is $G=\C^*$, and $\sigma(t)=\overline{t}^{-1}$. The fixed points of $\sigma$ is the circle $S^1$ of elements of norm 1 in $\C^*$. 

  If $Y$ is a $\C^*$-variety, a real form $\mu$ on $Y$ as a $\C^*$-variety which is compatible with $\sigma$ is called a  \emph{real circle form}, and the action of $(G,\sigma)$ on $(Y,\mu)$ is called a \emph{circle action on $Y$} or a \emph{real circle action on $Y$} . 

 Let $W_k$ be the two-dimensional $\C^*$-module with weights $(k,-k)$. In other words, $t(x,y)=(t^kx,t^{-k}y)$ for $t\in\C^*$.  
 
 Fix $m\ge 1$, and $n=2m+1$. 
 In this section, consider the following $\C^*$-variety :$Y=W_2\times W_n$.

The following  defines a real circle form on $Y$ :
$$ \mu_0(\left(\begin{array}{c}a \\b\end{array}\right), \left(\begin{array}{c}x \\y\end{array}\right))=(\left(\begin{array}{c}\overline{b} \\\overline{a}\end{array}\right),\left(\begin{array}{c}\overline{y} \\\overline{x}\end{array}\right)).$$

We call this form ``a linear circle action on affine space". Now we will construct a family of other real circle forms on $Y$. There might be many others, but we restrict the study to only some cases (where equivariant vector bundle methods can be used for constructing and distinguishing different forms). 

Later, we will explain how these real forms come from non-equivalent $O(2)$-actions found by G. Schwarz in \cite{S}.

Given $h\in\R[T]$, let  $M_h\in GL_2(\C[a,b])$ be the matrix given by
$$ M_h=\left(\begin{array}{cc}1-Th^2 & a^nh^n \\-b^nh^n & \sum_{j=0}^{n-1}(Th^2)^j\end{array}\right),
$$

Let $\varphi_h$ be  the automorphism of $Y$:
$$\varphi_h(\left(\begin{array}{c}a \\b\end{array}\right),\left(\begin{array}{c}x \\y\end{array}\right))=(\left(\begin{array}{c}a \\b\end{array}\right),M_h\left(\begin{array}{c}x \\y\end{array}\right)).$$
Note that $\varphi_h$ is $\C^*$-equivariant: $\varphi_h\in Aut_{\C^*}(Y)$.

We now state the main theorem of this article.

\begin{theorem} Let $h\in \R[T]$ be a real polynomial, and  let $\mu_h=\varphi_h\mu_0$.
\begin{itemize}
\item[(i)] $\mu_h$ defines a  real circle form on $Y=W_2\times W_{2m+1}$;
\item[(ii)] $\mu_h$ is equivalent to $\mu_{h'}$ if and only if there exists $r\in\R^*$ such that 
$$h(T)\equiv rh'(r^2T) \mod (T^m).$$
\end{itemize}
\end{theorem}

Note that, if one forgets the action, and one considers  $\mu_h$ as a real form of $Y$ as a variety, then they are all equivalent. 

\subsection{
Two special cases}

\begin{itemize}

\item[(1)] As a consequence, for $m=1$, we find two inequivalent real forms. For $h=0$, we find the linear circle action, and for $h=1$, we have another, non-equivalent form. This case is  precisely the action described in \cite{F-MJ}. 

\item[(2)] For $m=2$, if $h=c_0+c_1T$ and $h'=c_0'+c_1'T$, with all coefficients being real, the theorem implies that the real forms $\varphi_h\circ\mu_0$ and $\varphi_{h'}\circ\mu_0$ are equivalent if and only if : 
\begin{itemize}
\item[(i)] $c_0=c_0'=c_1=c_1'=0 $; 
\item[(ii)] $c_0=c_0'=0, c_1c_1'\not=0$; 
\item[(iii)] $c_1=c_1'=0, c_0c_0'\not=0$; or 
\item[(iv)] $c_0c_0'c_1c_1'\not=0$ and  $(c_0'/c_0)^3=c_1'/c_1$. 
\end{itemize}
In particular one finds an infinite family of inequivalent  real circle forms of the same $\C^*$-linear action.

\end{itemize}

\subsection{Three linear subgroups of $A=Aut_{\C^*} Y$}

In order to prove this theorem, we will use certain subgroups of $A=Aut_{\C^*} Y$. They are related to vector-bundle automorphisms described in general in the following sections.

(1) Consider the subgroup $\Lambda\subset GL_2(\C[a,b])$ defined as follows. Set $T=ab$. 
$$\Lambda=\{M=\left(\begin{array}{cc}P(T) & a^nQ(T) \\b^nS(T) & R(T)\end{array}\right)| P,Q,R,S\in\C[T], \det(M)=\Delta\in\C^*\}.$$

We define an action of $\Gamma$ on the group $\Lambda$. If 
$$M=\left(\begin{array}{cc}P(T) & a^nQ(T) \\b^nS(T) & R(T)\end{array}\right), \quad\text{then}\quad \ga M= \left(\begin{array}{cc}\overline{R}(T) & a^n\overline{S}(T) \\b^n\overline{Q}(T) & \overline{P}(T)\end{array}\right).$$

Now consider the subgroup $A_\Lambda\subset A$ of automorphisms  of the form 
$$\varphi_M((\left(\begin{array}{c}a \\b\end{array}\right),\left(\begin{array}{c}x \\y\end{array}\right))=(\left(\begin{array}{c}a \\b\end{array}\right),M\left(\begin{array}{c}x \\y\end{array}\right))$$
 where $M\in\Lambda$.

By construction, $\ga\varphi_M=\varphi_{\ga M}$, and $\varphi_M^{-1}=\varphi_{M^{-1}}$. 

(2) We consider also a subgroup of automorphisms of $Y$ coming from the linear circle action. More precisely, if $\omega\in\C$ is of norm 1, we construct the automorphism $\rho_\omega$ by 
$$\rho_\omega(\left(\begin{array}{c}a \\b\end{array}\right),\left(\begin{array}{c}x \\y\end{array}\right))=(\left(\begin{array}{c}\omega^2a \\\omega^{-2}b\end{array}\right),\left(\begin{array}{c}\omega^nx \\\omega^{-n}y\end{array}\right)).$$

(3) Finally, if $r\in \R$, we consider a subgroup of automorphisms that act  only on the first factor:
$\ell_r(\left(\begin{array}{c} a \\ b\end{array}\right),\left(\begin{array}{c}x \\ y\end{array}\right))=(\left(\begin{array}{c}ra \\rb\end{array}\right),\left(\begin{array}{c}x \\y\end{array}\right)).$

\bigskip\bigskip 
{\bf Proof the theorem, part (i)}:  It suffices to note that $\varphi_h\in A_\Lambda$, and that $M_h^{-1}=\ga M_h$. This implies that the cocycle condition $\ga\varphi_h^{-1}=\varphi_h$ is verified, and thererfore $\mu_h=\varphi_h\mu_0$ is an anti-holomorphic involution on $Y$, and it is compatible with $\sigma$, and therefore is a real circle form.\qed

\subsection{Equivalences}

The first step in proving the equivalence relation of part (ii) of the theorem is to reduce the problem of equivalence to a computation using automorphisms in $A_\Lambda$.

\begin{proposition} The real forms $\mu_h$ and $\mu_{h'}$ are equivalent if and only if there exists $\psi\in A_\Lambda$ and $r\in \R^*$ such that  $\psi\varphi_h\ga \psi^{-1}=\varphi_{h''}$ where $h''(T)=rh'(r^2T)$.
\end{proposition}

This result is a special case of the general Proposition 4.2, stated later. The proof for this case is given completely here.

\demo  First note that $\ell_r\varphi_{h'}\ga\ell_r^{-1}=\ell_r\varphi_{h'}\ell_r^{-1}=\varphi_{h''}$. Thus one direction of the equivalence is clear.

Now suppose that $\mu_h$ is equivalent to $\mu_{h'}$. This means that there is an automorphism $\psi\in A$ such that $\psi\circ\varphi_h\circ\ga \psi^{-1}=\varphi_{h'}$. 

We will now follow an argument from Masuda and Petrie \cite{MP2}   to create a vector bundle automorphism which has the same property.  Let $H=\{\pm 1\}$ be the subgroup of $\C^*$ of two elements. Since $\psi$ is equivariant, it stabilizes the fixed point set of $H$. That is, $\psi$ stabilizes the subvariety defined by $B=V(x,y)$. Also, by equivariance, and from the hypothesis that $\psi\circ\varphi_h\circ\ga \psi^{-1}=\varphi_{h'}$, we see that there exists $\lambda\in\C^*$ such that 
$$\psi(\left(\begin{array}{c}a \\b\end{array}\right),\left(\begin{array}{c}0\\0 \end{array}\right))=(\left(\begin{array}{c}\lambda a \\\overline{\lambda}b\end{array}\right),\left(\begin{array}{c}0 \\0\end{array}\right)). 
$$

By choosing an appropriate $r\in \R$ and $\omega\in\C$ of norm 1, one can construct $\psi_1=\ell_r^{-1}\rho_\omega\psi$ such that $\psi_1$ fixes pointwise $B$. Note that $\psi_1\circ\varphi_h\circ\ga\psi_1^{-1}=\varphi_{h''}$.

Finally $Y$ is the total space of the normal bundle of $B$ in $Y$. Thus $\psi_1$ induces an automorphism $\psi_2$ on the normal bundle, which gives therefore an automorphism in $A_\Lambda$. 

By the fact that $\varphi_h$ and $\varphi_{h'}$ are also in $A_\Lambda$, we have that $\psi_2\circ\varphi_h\ga\psi_2^{-1}=\varphi_{h''}$, and the proposition is proven.

\qed

Now we will show that the real circle actions $\mu_h$ which were constructed are all equivalent over a certain open set of $Y$. Then part (ii) of the theorem can be proven using this result and Proposition 3.2.

Notation: $V(ab)\subset Y$ is the closed subvariety $Y$ defined by the equation $ab=0$. 

Note first that $U=Y\setminus V(ab)$ is an open $G$-subvariety, and all the real forms $\mu_h$ restrict to real circle forms on $U$. We will start by showing that they are all equivalent on $U$.

Consider first the  analogous setting that was used on $Y$. Let $A'=Aut_{\C^*} U$.
Let  $\Lambda'\subset GL_2(\C[a,b,T^{-1}])$ be defined as follows. 
$$\Lambda'=\{M=\left(\begin{array}{cc}P(T) & a^nQ(T) \\b^nS(T) & R(T)\end{array}\right)| P,Q,R,S\in\C[T,T^{-1}], \det(M)=\Delta=cT^k, c\in\C^*,k\in\Z\}.$$

As before,  define an action of $\Gamma$ on the group $\Lambda'$. If 
$$M=\left(\begin{array}{cc}P(T) & a^nQ(T) \\b^nS(T) & R(T)\end{array}\right), \quad\text{then}\quad \ga M= \left(\begin{array}{cc}\overline{R}(T) & a^n\overline{S}(T) \\b^n\overline{Q}(T) & \overline{P}(T)\end{array}\right).$$

Consider the subgroup $A'_{\Lambda'}$ of $A'$ of automorphisms  of the form 
$$\varphi_M((\left(\begin{array}{c}a \\b\end{array}\right),\left(\begin{array}{c}x \\y\end{array}\right))=(\left(\begin{array}{c}a \\b\end{array}\right),M\left(\begin{array}{c}x \\y\end{array}\right))$$
 where $M\in\Lambda'$.

By construction, $\ga\varphi_M=\varphi_{\ga M}$, and $\varphi_M^{-1}=\varphi_{M^{-1}}$.

For $h\in \R[T]$, construct the matrix $K_h\in \Lambda'$ by
$$K_h=\left(\begin{array}{cc}1 & \frac{a^nh}{T^m} \\\frac{b^nh(\sum_{j=0}^{m-1}Th^2)}{T^m} & \sum_{j=0}^{m}Th^2\end{array}\right)
$$
and let $\psi_h'\in A'_{\Lambda'}$ be the corresponding automorphism of $U$:
$$ \psi_h'(\left(\begin{array}{c}a \\b\end{array}\right),\left(\begin{array}{c}x \\y\end{array}\right) )=(\left(\begin{array}{c}a \\b\end{array}\right),K_h\left(\begin{array}{c}x \\y\end{array}\right)).$$

\begin{lemma} Let $U=Y\setminus\{V(T)\}$. The real forms $\mu_h$ restricted to $U$ are all equivalent. More precisely, $\mu_h=\psi_h'\mu_0\psi_{h'}'^{-1}$. 
\end{lemma}

\demo It suffices to check that $\det(K_h)=1$ and that $K_h \ga K_h^{-1}=M_h$. 
\qed

\begin{lemma} Let $\Psi\in\Lambda'$ be a matrix such that $\ga\Psi=\Psi$ and such that $\det(\Psi)$ is a non-zero complex constant. Then $\Psi$ is of the form
$$\Psi=\left(\begin{array}{cc}\alpha & 0 \\0 & \overline{\alpha}\end{array}\right)$$
for some $\alpha\in\C^*$.
\end{lemma}

\demo Set $\Psi=\left(\begin{array}{cc}P(T) & a^nQ(T) \\b^nS(T) & R(T)\end{array}\right)$ where  $P,Q,R,S\in\C[T,T^{-1}]$, with  determinant $\Delta\in\C^*$. The first coefficient of $\Psi\ga\Psi^{-1}$ is $(P\overline{P}-T^nQ\overline{Q})/\overline{\Delta}$.  Therefore, we have the condition that $P\overline{P}-T^n Q\overline{Q}=\overline{\Delta}$. Consider $P\overline{P}$  and $Q\overline{Q}$ as Laurent polynomials. The highest monomial in $P\overline{P}$ is of even degree, and the highest monomial in $T^nQ\overline{Q}$ is of odd degree. Thus they do not cancel. The same holds for the lowest degree monomials. Thus the only possibility to have the given equality is that $Q=0$ and $P$ is of degree 0. Set $P=\alpha$. We find $\overline{\Delta}=\alpha\overline{\alpha}$, thus $\Delta$ is real. By considering the other coefficients of $\Psi$, we find similarly that $S=0$, and that $\Psi$ is of the given form.
\qed

\bigskip

{\bf Proof of Theorem 3.1}:  By Proposition 3.2, $\mu_h$ is equivalent to $\mu_{h'}$ if and only if there exists $r\in\R^*$ and $N\in\Lambda$ such that $N M_h \ga N^{-1}=M_{h''}$, where $h''(T)=rh'(r^2T)$. This means
$$N K_h \ga K_h^{-1} \ga N^{-1}=K_{h''}\ga K_{h''}^{-1}.$$

Set $\Phi=K_{h''}^{-1}NK_h\in\Lambda'$. $N$ satisfies the given equation if and only if $\ga\Phi=\Phi$. Also, since $\det(K_h)=\det(K_{h''})=1$, and since $\det{N}$ is a constant, the same holds of $\Phi$. By  Lemma 3.4, there exists $\alpha\in\C^*$ such that $\Phi$ is diagonal with terms $\alpha$ and $\overline{\alpha}$. 

The coefficients of N are all polynomials in $T$. Note that $K_h=\left(\begin{array}{cc}1 & a^nq_h \\ b^ns_h & r_h\end{array}\right)$ where $q_h=h/T^m$, $s_h=h(\sum_{j=0}^{m-1}Th^2))/T^m$, and $r_h=\sum_{j=0}^m (Th^2)^j$. The coefficients for $N$ are all polynomials in $\C[a,b]$ if and only if:
$$\alpha s_{h}r_{h''}-\overline{\alpha}r_hs_{h''}\in\C[T]$$ and

$$\alpha q_{h''}-\overline{\alpha}q_h\in\C[T].$$
First of all, if $h\equiv h'' (\mod T^m)$, then one can choose $\alpha=1$, and the two conditions are verified. That is, one can construct the matrix $N\in \Lambda$ which induces the equivalence of $\varphi_h$ and $\varphi_{h''}$. 

For the converse, note that 
 the second equation implies that $\alpha h \equiv\overline{\alpha}h'' (\mod T^m)$.

Since $h$ and $h''$ are real, if the real part of $\alpha$ is not zero, this implies that $h \equiv h'' (\mod T^m)$. If $\alpha$ is purely imaginary, then $h \equiv-h'' (\mod T^m).$ We also know that $\varphi_{h''}$ is equivalent to $\varphi_{-h''}$, by posing $r=-1$, which gives $rh''(r^2T)=-h''(T)$. Thus the existence of an appropriate matrix $N$ implies the condition in the theorem.

\qed

\section{Real forms of equivariant vector bundles}

The examples given here were constructed using techniques of equivariant $O_2(\C)$-vector bundles over $O_2(\C)$-modules. In \cite{K-S}, Kraft and Schwarz studied the following situation.   Let $G$ be a reductive complex group, and let $B$ and $F$ be two $G$-modules. The total space of an equivariant $G$-vector bundle over $B$ with zero fiber $F$, as a variety, is isomorphic to affine space. If an equivariant bundle is not trivial, it can happen that the associated action on the total space is non-linearizable. In fact, all examples of algebraic non-linearizable actions on affine complex space known at this time are constructed in this way. Kraft and Schwarz studied the case where the algebraic quotient of $B$ by $G$ is the affine line. They showed that the set of equivalence classes of $G$-vector bundles with base $B$ and zero fiber $F$ form a moduli space.  Another approach to the study of these $G$-vector bundles was developped in \cite{MP1}, \cite{MP2}, \cite{MMP}
 and \cite{MJ}, using the construction of equivariant bundles as subbundles of trivial ones. 
The method we used in Proposition 3.2 due to Masuda and Petrie was originally formulated  to show how non-equivalent $G$-vector bundles can lead to non-equivalent actions on the total spaces. 

In this section, we will compare the study of real forms of equivariant vector bundles and equivariant real forms on the total space. 
In the next section, we will describe how to use $O_2(\C)$-bundles to construct real forms of linear $\C^*$-actions on four-space.

Let $\pi: E\to B$ be a $G$-vector bundle. This means that it is a vector bundle over $B$, the group $G$ acts on $E$ and $B$, $\pi$ is equivariant, and $G$ acts linearly on the fibers. That is, for any $b\in B$ and $g\in G$, the restriction of the action of $g$ from $\pi^{-1}(b)$ to $\pi^{-1}(gb)$ is linear. If $f$ is a $G$-automorphism of $B$, a \emph {$G$-automorphism of $E$ over $f$}  is a $G$-automorphism $\varphi$ of $E$ such that $\varphi\circ\pi=\pi\circ f$, and such that $\varphi$ is linear on the fibers, that is, for any $b\in B$, the restriction of $\varphi$ from $\pi^{-1}(b)$ to $\pi^{-1}(f(b))$ is linear. 
 A \emph{ $G$-vector bundle automorphism} of $E\to B$ is a $G$-automorphism $\varphi$ of $E$ over the identity map on $B$.
  
 Suppose $B$ is a complex $G$-variety, $\sigma$ is a real form of $G$ and $\mu_B$ is a fixed real form of $B$ compatible with $\sigma$. 
 
 \begin{definition} A real form $\mu$ of $E$ is called a $G$-vector bundle real form compatible with $\sigma$ and $\mu_B$, if $\mu$ is an anti-holomorphic involution on $E$ defining a real form compatible with $\sigma$ and such that $\mu\circ\pi=\pi\circ\mu_B$, and $\mu$ is $\R$-linear on the fibers. Two $G$-vector bundle real forms  $\mu_1$ and $\mu_2$ are equivalent if there exists a $G$-vector bundle automorphism $\psi$ which conjugates $\mu_1$ to $\mu_2$. 
 \end{definition}   
 
 Clearly, $G$-vector bundle real forms compatible with $\sigma$ give real forms on the total space. Also, if two such real forms are equivalent as $G$-vector-bundle real forms, then they are equivalent as equivariant real forms on the total space. In general, the converse is not true. However, under some conditions, we can induce information of non-equivalence on the total space from non-equivalence as vector bundle real forms. The following proposition is essentially the result from Masuda-Petrie \cite{MP2} which was used in Proposition 3.2..

Suppose $\pi:E\to B$ is a trivial $G$-vector bundle and $f$ is a $G$-automorphism of $B$. Then $f$ induces a $G$-autormorphism on $E=B\times F$, which we also call $f$. Suppose now that $\sigma$ is a real form of $G$, and $\mu_B$ is a real form of $B$ compatible with $\sigma$, and $\mu$ is a $G$-vector bundle real form of $E$, compatible with $\sigma$ and $\mu_B$. If $f$ is a $G$-equivariant automorphism of $B$ which preserves $\mu_B$, (that is, $f\mu_B f^{-1}=\mu_B$), then $f^*\mu=f\mu f^{-1}$ is another vector bundle real form of $E$ compatible with $\sigma$ and $\mu_B$. By definition $f^*\mu$ and $\mu$ are equivalent as equivariant real forms of $E$.
 
 \begin{proposition} Let $G$ be an algebraic complex group with real form $\sigma$, and let $B$ be a $G$-variety with real form $\mu_B$ compatible with $\sigma$. Suppose  $\pi: E\to B$ is a trival $G$ vector bundle.
 
 Suppose that $G$ has a subgroup $H$ such that $E^H=B$, the zero section of the vector bundle.
 
  Then, given two equivariant vector bundle real forms $\mu_1$ anad $\mu_2$ on $E=B\times F$, they are equivalent as equivariant real forms of $E$ if and only if there exists a $G$-equivariant automorphism $f$ of $B$ which preserves $\mu_B$ such that $\mu_1$  and $f^*\mu_2$ are  equivalent as $G$-vector bundle real forms.
 \end{proposition}
 
 \demo Suppose $\psi$ is a $G$-equivariant automorphism of $E$ such that $\psi\circ\mu_1\circ \ga\psi^{-1}=\mu_2$. Since $\psi$ is equivariant for $H$, $\psi(B)=B$, where $B$ is identified with the zero section of the trivial vector bundle. $\psi_B=f$ is therefore a $G$-automorphism compatibile with $\sigma$ and $\mu_B$. Extending $f$ to $f=f\times id$ on $E=B\times F$, $(f^{-1})\circ\psi$  is a $G$-equivariant automorphism of $E$ which fixes pointwise the zero section $B$. This automorphism therefore induces a vector bundle automorphism $\psi'$ of the normal bundle of $B$ in the total space, which is $G$-isomorphic to $E$.  Since $\mu_1$ and $\mu_2$ are linear on the fibers, $\psi'\mu_1\ga{\psi'}^{-1}=f^*\mu_2$.
 \qed

\section{$O_2(\C)$-vector bundles}

In \cite{S} Schwarz gave a description of certain moduli spaces of $G$-vector bundle equivalence classes for $G=O_2(\C)=\C^*\rtimes\Z/2\Z$. The irreducible representations of $G$ are given by the one-dimensional ones and, for every $k\in\N\setminus\{0\}$, the two dimensional representation $W_k$ where $\C^*$ acts by weights $k$ and $-k$, and $G$ admits an involution which exchanges the two eigenspaces.  Schwarz,  showed in particular that the set of equivalence classes $VEC(W_2;W_{2m+1})$ of $G$-vector bundles with base $W_2$ and zero fiber $W_{2m+1}$ is a moduli space isomorphic to $\C^m$. All of these vector bundles come from considering the linear action of $\C^*$ on the total space, and then finding a compatible (regular) involution. One can then find non-linearizable $G$-actions on four-space.

In \cite{F-MJ} it was shown that, for the case $m=1$, one can use a certain involution from a non-linearizable action formed in this way  to construct an equivariant real circle form for the $\C^*$-action on the total space (see also \cite{W}). In order to do this, one needed to use that the involution commutes with conjugation. More precisely, the regular involution coming from a non-linearizable action of $G$ on affine four-space composed with conjugation defines an anti-holomorphic involution, since they commute. It was shown that the action given there was not equivalent to the linear circle action. 

In the present article, the cases treated come from families of non-equivalent vector bundles. To show that two real forms are not equivalent in general is more difficult than simply showing they are not equivalent to the linear one. The key to this calculation was to use that they are all equivalent on the open set $U$ constructed in section. 

We describe  now how to find examples of real circle actions on complex four-space with a linear action using $G$-vector bundles. Denote by $s\in G$ an order two element of $G$ which is not in the normal subgroup $\C^*$ : $G=\C^*\rtimes\{1,s\}$.  

A $G$-vector bundle $W_k\times W_n$ is constructed by defining the linear $\C^*$-action with weights $(k,-k,n,-n)$, and an action of the element $s$ respecting the conditions that $s^2=1$ and $st=t^{-1}s$. This restriction is very similar to that of finding real circle forms. 

For the trivial $G$-vector bundle $W_k\times W_n$  the action of  $s$ is given by the involution $\tau_0$ which exchanges the coordinates in the base and the zero fiber:

$$t(\left(\begin{array}{c}a \\b\end{array}\right),\left(\begin{array}{c}x \\y\end{array}\right))=(\left(\begin{array}{c}t^ka \\t^{-k}b\end{array}\right),\left(\begin{array}{c}t^nx \\t^{-n}y\end{array}\right))
$$
for $t\in\C^*$, and
$$
\tau_0(\left(\begin{array}{c}a \\b\end{array}\right)\left(\begin{array}{c}x \\y\end{array}\right) )=(\left(\begin{array}{c}b \\a\end{array}\right),\left(\begin{array}{c}y \\x\end{array}\right)).$$

All other $G$-vector bundle with base $W_k$ and zero fiber $W_n$ is equivalent to one where the $\C^*$-action remains the same, and the action of $s$ is given by an involution $\tau$ is of the form $\tau=\varphi\circ\tau_0$, where $\varphi$ is a $\C^*$-vector bundle automorphism, and $\tau$ is an involution. Analogously to the study of real forms, one can define an action of $\{1,s\}$ on the group of $\C^*$-vector bundle automorphisms on the trivial $\C^*$-vector bundle $W_k\times W_n$ by conjugation with $\tau_0$ :  $\sa\varphi=\tau_0\circ\varphi\circ\tau_0$. The condition for $\tau$ to define a $G$-vector bundle is then that $\sa\varphi^{-1}=\varphi$. Also, two $G$-vector bundles defined by the involution $\tau_1=\varphi_1\circ\tau_0$ and $\tau_2=\varphi_2\circ\tau_0$ are equivalent if and only if there exists a $\C^*$-vector bundle automorphism $\psi$ such that $\psi\varphi_1\sa\psi^{-1}=\varphi_2$. 

Note that the real form $\mu_0$ which defines a circle action on the total space is simply $\tau_0\circ conj= cong\circ\tau_0$, where $conj$ is the standard complex conjugation  on $\C^4$. Now suppose that one restricts  to vector bundle automorphisms $\varphi$ which commute with conjugation. Then $\mu=\varphi\circ\mu_0$ defines another $(\C^*,\sigma)$ equivariant real form. The question of equivalence, however, is not identical.  If $\psi$ can be constructed to give equivalence of $O_2(\C)$-vector bundles in such a way that it commutes with conjugation, then it gives an equivalence of the corresponding real circle forms. However, if $\psi$ does not commute with conjugation, the situation can become quite different. In the next section, such an example is treated. 

The explicit automorphisms used in this article were obtained by calculating the non-equivalent vector bundles using the descriptions from \cite{K-S} and from  \cite{MP1} . 

\section{ The $(1,2)$ case}

There are also non-trivial $O_2(\C)$-vector bundles whose base is $W_1$. These, however, do not in general lead to non-equivalent real circle forms of $\C^4$ with the corresponding $\C^*$ actions. First of all, the argument of Masuda and Petrie does not apply. But more than that, there are explicit  examples where an involution coming from a non-trivial equivariant vector bundle leads to a linearizable real circle action. (See also \cite{F-MJ}.)

We present here the following case. Consider the  $O_2(\C)$-vector bundle whose total space in $\C^4$ with weights $(1,-1,2,-2)$ for the $\C^*$-action.

Consider the involution $\tau$ defined by 
$$
\tau(\left(\begin{array}{c}a \\b\end{array}\right)\left(\begin{array}{c}x \\y\end{array}\right) )=\left(\begin{array}{c}b \\a\end{array}\right),\Phi\left(\begin{array}{c}y \\x\end{array}\right))$$

where 
$$\Phi=\left(\begin{array}{cc}1-T & a^4 \\-b^4 & 1+T+T^2+T^3\end{array}\right).$$

(Here, as before, $T=ab$.)  It is known that this involution together with the $\C^*$ linear action defines an $O_2(\C)$-vector bundle which is non-trivial. (It is not known whether the $O_2(\C)$-action on the total space is linearizable or not.)

Let $\tau_0$ be the involution which exchanges $a$ and $b$, and exchanges $x$ and $y$. Then $\tau=\varphi\circ\tau_0$, where $\varphi$ is a vector bundle involution: $\varphi(\left(\begin{array}{c}a \\b\end{array}\right),\left(\begin{array}{c}x \\y\end{array}\right))=(\left(\begin{array}{c}a \\b\end{array}\right),\Phi\left(\begin{array}{c}x \\y\end{array}\right) )$.

Now consider the total space $Y$ as a $\C^*$-variety. Since $\varphi$ commutes with complex conjugation, $\varphi\circ\mu_0=\mu$ defines a circle form on $Y$. 

\begin{proposition} The circle form $\mu=\varphi\mu_0$ on $Y=W_1\times W_2$ is equivalent to the linear circle action $\mu_0$. 
\end{proposition}

\demo

Let $N\in GL_2(\C[a,b])$ be the matrix:
$$N=\left(\begin{array}{cc}1-\frac{(1-i)}{2}T-\frac{(1+i)}{4}T^2 & \frac{a^4}{4}(1-i) \\\frac{-b^4}{4}(3-i+(1+i)T) & 1+\frac{1-i}{4}(2T+T^2+T^3)\end{array}\right).
$$
Then 
$$\psi(\left(\begin{array}{c}a \\b\end{array}\right),\left(\begin{array}{c}x \\y\end{array}\right) )=(\left(\begin{array}{c}a \\b\end{array}\right), N\left(\begin{array}{c}x \\y\end{array}\right))$$
is a $\C^*$-automorphism.  

Also $\psi\mu_0\psi^{-1}=\mu$, since $N\ga N^{-1}=\Phi$.  Therefore, $\mu_0$ is equivalent to $\mu$.

\qed

Note that the coefficients of $N$ are not real, that is, $\psi$ does not commute with complex conjugation.

\section{Questions}

\subsection{Quotients}

An interesting question to consider is the real forms of the quotients. Given $Y=W_2\times W_n$, the algebraic quotient $Z=Y//(\C^*)$ is a singular threefold. Set $T=ab$, $W=xy$, $U=a^ny^2$ and $V=b^nx^2$. These polynomials generate the invariant ring, and they satisfiy the equation $UV-T^nW^2=0$. More precisely, $Z\cong V(UV-T^nW^2)\subset\C^4$. Any  circle form on $Y$ induces a real form on $Z$. Thus it is natural to ask if inequivalent circle forms on $Y$ induce inequivalent real forms on $Z$. The answer is not obvious, because one would have to show that the automorphism group of $Z$ has restrictions, allowing to lift them to equivariant automorphisms on $Y$. In particular, it is important to know if  automorphisms of $Z$ preserve the stratification  given by the quotient map $Y\to Z$. 

Note also, that one can show that  the real forms are all equivalent in the diffeomorphic category. 

\subsection{Other examples}

In this article, all the real forms that are considered come from a very particular type. They are constructed from regular involutions of non-trivial equivariant vector bundles, and these involutions all commute with complex conjugation  $\C^4$ (as a variety). 

It would be interesting to consider other cases. In particular, it is possible that, in comparison to the case presented in Section 6, there are examples of equivalent equivariant $O_2(\C)$-bundles that lead to inequivalent real circle forms. The techniques applied using  $O_2(\C)$-vector bundles would have to be adapted to study this case.

\bigskip
\bigskip 
\noindent L. Moser-Jauslin\\
Universit\'e de Bourgogne Franche-Comt\'e \\
Institut de Math\'ematiques de Bourgogne -UMR 5584 du CNRS\\
 9 av. Alain Savary, BP 47870 \\ Dijon 21078, France\\
e-mail : moser@u-bourgogne.fr

\end{document}